\documentclass{amsart}

\usepackage[centertags]{amsmath}
\usepackage{amssymb}
\usepackage{amsfonts}

\newtheorem{theorem}{Theorem}
\newtheorem{lemma}[theorem]{Lemma}
\newtheorem{propo}[theorem]{Proposition}

\newtheorem{rmk}[theorem]{Remark}

\catcode`\@=11
\def\theequation{\@arabic{\c@equation}}
\def\thetheorem{\@arabic{\c@theorem}}

\newcommand{\bt}{\begin{theorem}}
\newcommand{\et}{\end{theorem}}
\newcommand{\bl}{\begin{lemma}}
\newcommand{\el}{\end{lemma}}
\newcommand{\bp}{\begin{propo}}
\newcommand{\ep}{\end{propo}}
\newcommand{\bc}{\begin{coro}}
\newcommand{\ec}{\end{coro}}
\newcommand{\bdefn}{\begin{defn}}
\newcommand{\edefn}{\end{defn}}
\newcommand{\bcas}{\begin{cases}}
\newcommand{\ecas}{\end{cases}}
\newcommand{\br}{\begin{rmk}}
\newcommand{\er}{\end{rmk}}
\newcommand{\be}{\begin{equation}}
\newcommand{\ee}{\end{equation}}
\newcommand{\ben}{\begin{equation*}}
\newcommand{\een}{\end{equation*}}
\newcommand{\ba}{\begin{array}}
\newcommand{\ea}{\end{array}}
\newcommand{\bea}{\begin{eqnarray}}
\newcommand{\eea}{\end{eqnarray}}
\newcommand{\bean}{\begin{eqnarray*}}
\newcommand{\eean}{\end{eqnarray*}}
\newcommand{\beit}{\begin{itemize}}
\newcommand{\eeit}{\end{itemize}}
\newcommand{\benu}{\begin{enumerate}}
\newcommand{\eenu}{\end{enumerate}}

\newcommand{\ds}{\displaystyle}

\newcommand{\e}{\varepsilon}
\def \R {\mathbb{R}}
\def \N {\mathbb{N}}

\newcommand{\into}{\int_\Omega}

\begin{document}

\title{A note on a non-linear Krein-Rutman theorem}

\author{Rajesh Mahadevan}
\address{Departamento de Matem\'atica, Universidad de Concepci\'on,
  Casilla 160-C, Concepci\'on, CHILE}
\email{rmahadevan@udec.cl}
\thanks{This work was funded partly by CMM, Universidad de Chile 
while the author was a post-doc and partly, by the University of
Concepci\'on where the author is now a faculty member of the
Department of Mathematics.}


\date{September 28, 2006.}


\begin{abstract}
In this note we will present an extension of the Krein-Rutman theorem
for an abstract non-linear, compact, positively  1-homogeneous
operators on a Banach space having the properties of being increasing
with respect to a convex cone $K$ and such that there is a non-zero $u
\in K$ for which $M \, Tu \succcurlyeq u$ for some positive constant
$M$. This will provide a uniform framework for recovering the
Krein-Rutman-like theorems proved for many non-linear differential
operators of elliptic type, like the $p$-Laplacian cf. Anane \cite{A},
Hardy-Sobolev operator cf. Sreenadh \cite{Sr}, Pucci's operator
cf. Felmer et. al. \cite{FQ}. Our proof follows the same lines as in
the linear case cf. Rabinowitz \cite{R} and is based on a bifurcation
theorem.  

\end{abstract}

\maketitle

\section{Introduction}

\noindent
Let $X$ be a real Banach space. Let $K$ be a closed convex cone in $X$
with vertex at $0$, that is a set  having  the properties: 

\noindent
(i) $0 \in K$, (ii) $x \in K, t \in \R^+ \implies t \ x \in K$, and
(iii) $x,y \in K \implies x + y \in K$. 

\noindent
We further assume that 
\begin{equation*}
\hspace{-5cm} (A) \hspace{5cm} \quad K \cap {-K}   = \{0\}. 
\end{equation*}  
The cone $K$ induces an {\em ordering} $\preccurlyeq$ on $X$ defined
as follows. Given any $x,y \in X$ we say that
\begin{equation}\label{ord}
x \preccurlyeq y \Leftrightarrow  y-x \in K.
\end{equation}
The ordering in (\ref{ord}) is said to be {\em strict} if $x
\preccurlyeq y$ and $x \neq y$ and this will be denoted by $x \prec
y$. A mapping $T: X \to X$ is said to be {\em   increasing } if
$x \preccurlyeq y \implies T(x) \preccurlyeq  T(y)$ and it is said to
be {\em strictly increasing} if $x \prec y$ implies $Tx \prec Ty$. The
mapping is said to be {\em compact  } if it takes bounded subsets
of $X$ into relatively compact subsets of $X$. We say that the
mapping is {\em positively  1-homogeneous} if it satisfies the 
relation $T(t \, x) = t\, T(x) $ for all $x \in X$ and $t \in
\R^+$. We say that a real number $\lambda$ is an eigenvalue of the
operator if there exists a non-zero $x \in X$ such that $\lambda Tx =
x$. 

\medskip 
\noindent
Let us end this section by stating a simple and obvious fact concerning
closed convex sets which will be used later on.
\bl \label{ml} Given any $0 \prec x$ and $y \notin K$ there
exists a unique $\delta$ with $0 \le \delta(y)$ such that  
\bean
x + \lambda \, y & \in & K \mbox{ if } 0 \le \lambda \le \delta(y)\\
x + \lambda \, y & \notin & K \mbox{ if } \lambda > \delta(y).
\eean
Furthermore, we shall have $\delta >0$ if $x$ belongs to the interior
of $K$ denoted by $\mbox{\r{K}}$.
\el

\section{Main Theorem}
\setcounter{equation}{0}
\noindent
Let $X$ be a real Banach space and let $K$ be a closed convex cone
satisfying the assumption (A) given in the previous section.   
\bt \label{mt} Let $T: X \to X$ be an increasing, positively
$1$-homogeneous compact continuous operator(non-linear) on $X$ for
which there exists a non-zero $u \in K$ and $M>0$ such that 
\ben
\hspace{-5cm} (H) \hspace{5cm} M \,Tu \succcurlyeq u. 
\een
Then, $T$ has a non-zero eigenvector $x_0$ in $\mbox{K}$. Furthermore,
if $K$ has non-empty interior and if $T$ maps $K \setminus\{0\}$  into
$\mbox{\r{K}}$ and is strictly increasing, then $x_0$ is the unique
positive eigenvector in $K$ upto a multiplicative constant. And,
finally if $\mu_0$ be the corresponding eigenvalue, then it can be
characterized as the eigenvalue having the smallest absolute value and
furthermore, it is simple. 
\et
\noindent
We prove this theorem in exactly the same way as it is done in the
notes of Rabinowitz \cite{R} for linear operators and relies on the
following result cf. Rabinowitz \cite{R}, Corollaire 1. 
\bp  Let $X$ and $K$ be as in the statement of the previous
theorem. Let us be given a mapping $F: \R \times K \to K$ which is
compact, continuous and such that $F(0,x)=0$ for all $x \in K$. Then,
the equation $F(\lambda,x)=x$ has a non-trivial connected unbounded
component of solutions $\mathcal{C}^+$ in $\R^+ \times K$ containing
the point $(0,0)$.  
\ep
\medskip
\noindent
{\bf Proof of Theorem \ref{mt}: }First we prove the {\em existence} of
a positive eigenvector. 

\noindent
\underline{Step 1 : } Let $u \in K$ ($u \neq 0)$ be as in the
hypothesis (H) of the theorem. Let $\e >0$ be a parameter. We define a
parametrized family of operators $F_\e: \R^+ \times K \to K$ as follows
\be \label{parfam}
\ds{F_\e(\lambda,x) : =\lambda \, T \left( x+ \e \, u\right).}
\ee
Then, since $T$ is compact and continuous, each of these operators
$F_\e$ is clearly compact and continuous on $\R^+ \times
K$. Also they map $\R^+ \times K$ into $K$ since $K+\e \,u \subset K$
and $T$ maps $K$ into itself, which follows from the fact that $T$ is
increasing and $T(0)=0$. Let, by the above proposition,
$\mathcal{C}^+_\e \subset \R^+ \times K$ be a connected unbounded
branch of solutions to the equation   
\be\label{mt-ref1}
F_\e(\lambda,x) (:=\lambda T(x+ \e \, u)) = x.
\ee
{\bf Claim: } Then $\mathcal{C}^+_\e \subset [0, M] \times K$ for all
$\e >0$.

\smallskip
\noindent
We shall prove the claim by proving the following estimate: whenever,
$(\lambda,x) \in \mathcal{C}^+_\e$ we have
\be\label{mt-ref2}
\ds{x \succcurlyeq \left(\frac{\lambda}{M}\right)^n \ \e u \quad
  \forall n   \in \N.}
\ee  
Indeed, if the estimate (\ref{mt-ref2}) is true then we cannot have
$\lambda > M$. Otherwise, letting $n \to +\infty$ we will have $u \in
- K$ and hence, $u \in K \cap (-K)$ contradicting the hypothesis
(A). So, we are left to prove the estimate (\ref{mt-ref2}). We shall
do it by induction. Let $(\lambda,x)$ belong to
$\mathcal{C}^+_\e$. Then, $x = \lambda T(x +\e u)$ and we obtain, from
the properties of $T$ and the inequalities $x+ \e \, u \succcurlyeq \e
\, u$,  $x+\e \, u \succcurlyeq x$, respectively that
\bea
\label{mt-ineq1} x & \succcurlyeq & \lambda \, \e \, Tu\,, \\ 
\label{mt-ineq2} x & \succcurlyeq & \lambda Tx
\eea 
Now, using (H), it follows from (\ref{mt-ineq1}) that 
\be\label{mt-iter1}
\ds{x \succcurlyeq \lambda \, \frac{\e}{M} \, u}.
\ee
and we obtain our claim for $n=1$. Let us now assume that
(\ref{mt-ref2}) is true for $n=m$. Operating $T$ on
(\ref{mt-ref2}), we obtain 
\ben 
\ds{Tx \succcurlyeq \left(\frac{\lambda}{M}\right)^m \, \e Tu
  \succcurlyeq \left(\frac{\lambda^m}{M^{m+1}}\right) \, \e \, u} 
\een
So, using (\ref{mt-ineq2}), we obtain 
\ben
\ds{x \succcurlyeq \left(\frac{\lambda}{M} \right)^{m+1}\, \e \, u. }
\een
This, completes the induction step and proves (\ref{mt-ref2}).

\noindent
\underline{Step 2: }As we have shown that $\mathcal{C}^+_\e \subset 
[0,M] \times K$ for every $\e >0$ and besides, the branch is connected
and infinite starting from $(0,0)$, there must necessarily exist $x_\e$
with $\|x_\e\|=1$ and $\lambda_\e \in [0,M]$ such that
$(\lambda_\e,x_\e) \in \mathcal{C}^+_\e$. That is we have,
\be\label{mt-ref5}
x_\e=\lambda_\e \, T(x_\e + \e \, u), \quad \|x_\e\| = 1.
\ee 
Now, if we consider a sequence $\e$ which tends to $0$, then the
sequence $x_\e +\e u$ is bounded in $X$. As the operator $T$ is
compact, we may assume that $T(x_\e + \e u)$ converges for a 
subsequence and we suppose also that $\lambda_\e$ converges to some
$\lambda_0$ for an induced subsequence. Without loss of generality,
indexing this subsequence again by the same $\e$, we obtain that
$x_\e$ converges to some $x_0 \in K$ with $\|x_0\|=1$. Further, since
$T$ is a continuous operator, letting $\e \to 0$ in (\ref{mt-ref5}) we
obtain 
\be\label{mt-evp}
x_0=\lambda_0 \, Tx_0\,, \quad \|x_0\| = 1\,, x_0 \in K.
\ee
Thus, we have shown the existence of an eigenvector in
$K\setminus\{0\}$. Further, it follows from (\ref{mt-evp}) that
$\lambda_0 \neq 0$. 

\medskip
\noindent
From now on we shall assume that $\mbox{\em \r{K}}$ is non-empty and
that $T$ maps $K \setminus \{0\}$ into $\mbox{\em \r{K}}$ and is
strictly increasing. 

\medskip
\noindent
\underline{Step 3: }We shall now prove that $x_0$ obtained above is
the unique eigenvector in $K$ upto a multiplicative constant. Since,
$x_0$ is non-zero and we have assumed that $T$ maps $K \setminus
\{0\}$ into $\mbox{\em \r{K}}$ we obtain first of all that $x_0 \in
\mbox{\em \r{K}}$. By the same arguments, if $y \in K$ is  any other
eigenvector having eigenvalue $\lambda$ then, $\lambda \neq 0$ and $y
\in \mbox{\em \r{K}}$. First, using the fact that $x \in
\mbox{\em \r{K}}$ and $- y \notin K$, we apply Lemma \ref{ml}
and obtian $\delta(-y) > 0$ such that $x_0 - \delta(-y) \, y \in
K$. We claim that $x_0 = \delta(-y) \, y$. Otherwise, since $T$ is
strictly increasing, we will have  
\ben
T x_0 \succ T(\delta(-y) \, y).
\een
From this we deduce that $(\lambda_0)^{-1} x_0 \succ (\lambda)^{-1}
\delta(-y) y$ which, again by Lemma \ref{ml}, yields 
\be \label{mt-ref6}
\lambda_0 < \lambda.
\ee  
On the other hand, similarly, starting from the fact that $y -
\delta(-x_0) x_0 \in K$, we obtian $(\lambda)^{-1} \, y
\succcurlyeq (\lambda_0)^{-1} \, \delta(-x_0)$ yielding 
\be \label{mt-ref7}
\lambda \le \lambda_0.
\ee
Since (\ref{mt-ref6}) and (\ref{mt-ref7}) contradict each other, it
follows that our claim $x_0 = \delta(-y) \, y$  must be true. This
shows that any other eigenvector in $K$ is a multiple of $x_0$.

\medskip
\noindent
\underline{Step 4: }Let us now show the simplicity of
$\lambda_0$. From the arguments of the previous paragraph if $y$ is
any other eigenvector corresponding to $\lambda_0$ and if $y \in K$
then $y$ is a multiple of $x_0$. So, we deal now with the case that $y
\notin K$. We claim that $x_0 = \delta(y) \, y$. If not, arguing as
in the previous paragraph, we will obtain $\lambda_0 < \lambda_0$
which is absurd. Thus, we have proved that $\lambda_0$ is simple.

\medskip 
\noindent 
\underline{Step 5: }Let us now show that $\lambda_0$ has the smallest
absolute value among all eigenvalues of the operator $T$. Let $\lambda
\in \R$, different from $\lambda_0$, be any eigenvalue of $T$ with
corresponding eigenvector $y$. Since, $\lambda \neq \lambda_0$, this
would mean that neither $y$ nor $-y$ belongs to $K$. So, we have, $x_0
\pm \delta(\pm y ) \, y \in K$, which gives 
\be\label{mt-ref8}
\ds{x_0 \pm \frac{\lambda_0}{\lambda} \delta(\pm y) y \in K.}
\ee
Let us first consider the case $\lambda >0$. From (\ref{mt-ref8}),
using Lemma \ref{ml}, we get $\lambda_0 \le \lambda$ showing that
$\lambda_0$ has a smaller absolute value. In the case when $\lambda <
0$, from (\ref{mt-ref8}), we can obtain the inequalities 
\be\label{mt-ref9}
\ds{\frac{\lambda_0}{-\lambda} \, \delta(+ y ) \le \delta(- y) \,, 
\quad \mbox{ and }  \frac{\lambda_0}{- \lambda} \, \delta(- y ) \le 
\delta(+ y).}
\ee
From the above, we get 
\ben
\ds{\frac{\lambda_0}{-\lambda} \le \frac{\delta(- y)}{\delta(+ y)} \le
  \frac{- \lambda}{\lambda_0}}.   
\een
This, gives $|\lambda_0| \le |\lambda|$. Thus, we have shown that
$\lambda_0$ has the smallest absolute value among all eigenvalues of
$T$. \qed

\noindent
\br If the operator $T$ maps $K\setminus\{0\}$ into \r{K} then it is
called {\em strongly positive}. For such operators, the hypothesis (H)
holds. This can be shown easily by contradiction. However, it is not
necessary that $T$ be strongly positive for the hypothesis (H) to hold.
In the next section, we shall see some examples of operators which
fail to be strongly positive in the first place but we shall obtain
the existence of positive eigenvectors as they satisfy the hypothesis
(H). As far as the uniqueness and simplicity is concerned, we still
need the operator to be strongly positive and strictly increasing. For
a linear operator, strong positivity implies, automatically, that it
is strictly increasing. \qed     
\er

\section{Applications}
\setcounter{equation}{0}
\noindent
It has been shown in several papers \cite{A,FQ,Sr} that many
non-linear elliptic operators such as the $p$-Laplacian operator,
Hardy-Sobolev operator, Pucci's maximal operators etc. verify the
Krein-Rutman theorem. The methods that have been used to prove a
Krein-Rutman theorem for the $p$-Laplace operator or the Hardy-Sobolev
operator have relied very much on the variational structure of the
operators and on the use of special identities such as the Picone's
identity \cite{AH}, thus marking a complete difference with the
methods used for proving the same result for the Pucci's maximal
operators which are completely non-linear and have no variational
structure. It is in this context that our main theorem comes to
illustrate the main features that a non-linear operator should
possess in order to have a unique(upto a multiplicative constant)
positive eigenvector. 

The {\em existence part} is guaranteed as long as the operator is
positively $1$-homogeneous, is monotone with respect to a convex cone
which does  not necessarily have non-empty interior, and satisfies the
condition (H) of the theorem. When it comes to checking whether an
operator satisfies the hypotheses of the theorem, it is quite clear
that homogeneity is a straight forward condition to verify while,
the monotonicity depends on the operator space and the cone
chosen. Often, the natural choice is the cone of non-negative
functions in the function space and verifying this condition is,
in general, a task of seeing whether the weak comparison principle
holds. This is, usually, not a difficult problem. Verifying the
condition (H) on the other hand may involve a little more work and for
this, although not necessary, it is helpful if the operator satisfies
a strong maximum principle.   

The {\em uniqueness part} requires the operator to be strictly
increasing which is easily satisfied when the operator is increasing and 
injective. It also requires that the operator be strongly
positive. This is a tougher condition to verfiy and to begin with, the
cone definition has to be changed. A cone which has non-empy interior
in the space of continuous functions on a domain $\Omega$ vanishing on the
boundary is, for example, 
\be \label{regcon}
\ds{K^*=\left\{w \in C^1(\overline{\Omega})\,,  w \ge 0 
\mbox{ in } \Omega\,, \frac{\partial w}{\partial n} \le 0 \mbox{ on }  
\partial \Omega \right\}}.
\ee
To be able to work this cone the operator should have enough
regularity. Admitting this, the operator shall be strongly positive 
for this cone if, for example, it satisfies a Hopf maximum principle
which readily follows from a strong maximum principle if the boundary
of the domain satisfies an interior sphere condition. Thus, to sum up,
our conclusion is that the full Krein-Rutman theorem  
must be available for many non-linear elliptic operators defined on
fairly regular domains satisfying the weak comparison principle,
having good regularity properties and a strong maximum principle. We
shall illustrate this in the case of the $p$-Laplace operator,
Hardy-Sobolev operator and the Pucci operators.

\medskip
\noindent
Let $\Omega$ be a bounded domain having a smooth, connected boundary
in $\R^n$. Let us look at the Dirichlet eigenvalue problem for the
$p$-Laplace operator, $- \Delta_p \cdot:=\mathrm{div}(|\nabla \cdot
|^{p-2} \nabla \cdot) $ for $1 < p < \infty$ :
\be\label{evp-plap}
\left\{\ba{rcll}
- \Delta_p u & = & \lambda \, |u|^{p-2} u & \mbox{ in } \Omega \\
u & = & 0  & \mbox{ on } \Omega\,.
\ea\right.
\ee 
Or, now, let $0$ belong to the domain and let us look at the Dirichlet
eigenvalue problem for more general operators called the Hardy-Sobolev
operators. For this recall the Hardy-Sobolev inequality  \cite{H, GP} 
\be\label{HSineq}
\into \, |\nabla u|^p  dx \ge \ds{\left(\frac{n-p}{p}\right)^p} \into \,
w(x) |u|^p dx \mbox{ for all } u \in W^{1,p}_0(\Omega).
\ee
where $w$ is the weight function 
$$
w=\left\{ 
\ba{l}
\ds{\frac{1}{\left(|x|\log\frac{R}{|x|}\right)^n}} \mbox{ when } p=n\\
\\
 \ds{\frac{1}{|x|^p}} \mbox{ when }1<p<n\,,\\
\ea\right.
$$
and
$\ds{\left(\frac{n-p}{p}\right)^p}$ is the best constant for the
inequality. Let $\mu$ be a number strictly smaller than the best
constant. Then one defines the Hardy-Sobolev operator $L_\mu$ by
setting $L_\mu := -\Delta_p u - \mu w(x) |u|^{p-2} u$. In particular,
we recover the $p$-Laplacian operator if we take $\mu=0$. We may
consider the Dirchlet eigenvalue problem 
\be\label{HSevp}
L_\mu u = \lambda V(x) \, |u|^{p-2}u \mbox{ in } \Omega\,, u \in
W^{1,p}_0(\Omega)\,,u \neq 0.   
\ee
for a positive singular weight function $V$ whose singularity at $0$ 
is not worse than that of $w$ above, in the sense that, $\ds{\lim_{x \to 0}
w(x) V(x) = 0}$, as in \cite{NM,Sr}. The full Krein-Rutman for this 
problem has been proved, using variational techniques, by Sreenadh
\cite{Sr}. We now show that the same can be recovered from our theorem
by introducing a suitable operator framework. For this let us consider
the weighted reflexive Banach space $L^p(\Omega, V)$ and define  a
non-linear operator $T$ on it by setting $\ds{Tf := (L_\mu)^{-1}
  \left( V(x) |f|^{p-2} f \right) }$. It can be shown that this is a
well defined, compact, continuous operator by using arguments  similar
to those used in \cite{NM}. In fact, it has been shown in \cite{NM}
that if $f \in W^{1,p}_0(\Omega)$, then $V(x) |f|^{p-2} f \in
W^{-1,q}(\Omega)$ for the above class of singular weights $V$. In
order that $(L_\mu)^{-1}$ be well defined we need to examine the
existence and uniqueness of solution for the Dirichlet problem  
\be\label{HS-Dirpb}
L_\mu u = g \mbox{ in } \Omega\,, u \in W^{1,p}_0(\Omega)   
\ee
given any $g \in  W^{-1,q}(\Omega)$. A solution to this problem can be
obtained  by minimizing the energy functional 
$$
J(v)= \ds{\frac{1}{p} \int_\Omega |\nabla v |^p dx 
- \mu \frac{1}{p} \int_\Omega \frac{|v |^p}{|x|^p} dx - \int_\Omega g
v dx} \, . 
$$
By the Hardy-Sobolev inequality (\ref{HSineq}) and the fact that $\mu$
is smaller than the best constant in (\ref{HSineq}), it follows that the
functional is coercive. However, the problem is that it is not weakly
lower semicontinuous on $W^{1,p}_0(\Omega)$ as, except when $p=2$, the
functional is not convex for any $\mu > 0$.  Thus to show that the
minimum is attained one has to exhibit a minimizing sequence which
converges strongly in $W^{1,p}_0(\Omega)$. For this, one chooses,
using the Ekeland variational principle, a minimizing sequence $v_n$
for which $J'(v_n) \to 0$ and it can be checked, similarly as in
\cite{NM}, that such a sequence $v_n$ converges strongly in
$W^{1,p}_0(\Omega)$. Also, it can be argued  that what is important
for the theorem is that $L_\mu^{-1}$ be well defined on the positive
cone and and the desired uniqueness of positive solutions follows from
the weak comparison principle (WCP) proved for positive solutions of
such operators in \cite{CT}. The same WCP also shows that the operator
$T$ is monotone increasing on the positive cone. Now, it is quite easy
to check that $T$ is positively  $1$-homogeneous and its compactness
follows from the fact the image of this operator is contained in
$W^{1,p}_0(\Omega)$ which is compactly imbedded in $L^p(\Omega, V)$  
(cf. \cite{NM}). Therefore, it remains to verify the condition (H) to
be able to obtain the existence of a positive eigenvector. For this,
we use the $C^{1,\alpha}$ regularity of solutions to the Dirichlet
problem (\ref{HS-Dirpb}) proved by Tolksdorf \cite{T2} and the strong
maximum principle for such an operator \cite{T1,V}.  In fact, if $f$
is any non-negative smooth function with compact support then, by the
above two results, $Tf$ will be a strictly positive function on
$\Omega$ and we obtain (H) for a suitable constant $M$. If now we
would like to obtain also the {\em uniqueness} (upto a multiplicative
constant) of positive eigenvector, we need to change the setting to
the space of continuous functions vanishing on the boundary. Using the
$C^{1,\alpha}$ regularity and the Arzela-Ascoli theorem, we obtain the
compactness of the operator $T$.  The operator $T$ is strictly
increasing because it is increasing and injective. Finally, $T$ sends
$K^*\setminus\{0\}$ into the interior of $K^*$ because of the strong
maximum principle. So, the full Krein-Rutman theorem follows.

\medskip
\noindent
Now we apply the results of our theorem to the case of fully
non-linear elliptic operators where no variational methods can be
used. Consider the Dirichlet eigenvalue problem 
\be\label{fnlevp}
\ba{lcll}
F(D^2 u) & = & \lambda u & \mbox{ in } \Omega\\
u & = & 0 & \mbox{ on } \partial \Omega\,
\ea
\ee 
where $F$ is a fully non-linear elliptic operator like in the monograph
of Caffarelli and Cabr\'e \cite{CC} and which we assume, further, is
positively $1$-homogeneous. For example, let us consider the Pucci's
maximal operators $\mathcal{M}^{\pm}_{\lambda,\Lambda}$. The Pucci's 
maximal operator $\mathcal{M}^+_{\lambda,\Lambda}$
($\mathcal{M}^-_{\lambda,\Lambda}$)  is convex (respectively, concave) 
and hence, viscosity solutions of 
\be\label{fnl-Dirpb}
\ba{lcll}
F(D^2 u) & = & f & \mbox{ in } \Omega\\
u & = & 0 & \mbox{ on } \partial \Omega\,
\ea
\ee 
are in fact $C^{2,\alpha}$ regular  (cf. Section 6 \cite{CC}). So, we
may look at the eigenvalue problem as corresponding to the eigenvalue
problem $\lambda T u = u$ for the solution operator of the 
Dirichlet problem (\ref{fnl-Dirpb}) in $C^2(\Omega)$. The compactness
of $T$ follows from the $C^{2,\alpha}$ regularity. Also, it is known
that the weak comparison  principle holds for these operators
(cf. Proposition 2.2. \cite{FQ}). The strong maximum principle also 
holds (cf. Proposition 4.9 \cite{CC}). For all these reasons, we may
now apply our theorem and conclude the existence of a unique (upto a
multiplicative constant) positive eigenvector for these operators.

\medskip
\noindent
{\bf Acknowledgements: } The author would like to thank P. Felmer and
J. Davila of DIM, U. Chile and M. Montenegro of U. Campinas, Brazil
for many fruitful discussions. The author would also to thank the
anonymous referee for  valuable suggestions.

\bibliographystyle{amsplain}

\end{document}